\documentclass[conference]{IEEEtran}

\usepackage{fontawesome}
\usepackage{hyperref}

\def\pscircled#1{\textcircled{\resizebox{.5em}{!}{#1}}}

\IEEEoverridecommandlockouts
\usepackage{cite}
\usepackage{amsmath,amssymb,amsfonts}
\usepackage{amsthm}
\usepackage{algorithmic}
\usepackage{graphicx}
\usepackage{textcomp}
\usepackage{xcolor}
\newtheorem{theorem}{Theorem}
\newtheorem{lemma}[theorem]{Lemma}
\newtheorem{definition}{Definition}
\def\BibTeX{{\rm B\kern-.05em{\sc i\kern-.025em b}\kern-.08em
    T\kern-.1667em\lower.7ex\hbox{E}\kern-.125emX}}

\begin{document}

\title{Controlled Optimization with a Prescribed Finite-Time Convergence Using a Time Varying Feedback Gradient Flow\\
\thanks{ \small OA and YC were supported by the Center for Methane Emission Research and Innovation (\href{http://methane.ucmerced.edu}{\tt CMERI}) through the Climate Action Seed Funds grant (2023-2026) at University of California, Merced. 
\pscircled{\faDesktop}~\href{http://mechatronics.ucmerced.edu}{\tt mechatronics.ucmerced.edu}\par
NSO is supported by  The Scientific and Technical Research Councilof TÜRKİYE (TUBITAK) 2219-International Postdoctoral Research Fellowship Program for Turkish Citizens entitled ``Applied Fractional Calculus for Complex Big Data Monitoring and Control".  }
}

\author{Osama F. Abdel Aal$^{1}$, Necdet Sinan Özbek$^{1,2}$, Jairo Viola$^{1}$ and YangQuan Chen$^{1}$
\thanks{ \small
$^{1}$Dept. of Mechanical Engineering, University of California, Merced,
        5200 N. Lake Rd, Merced, CA 95343, USA
        \pscircled{\faEnvelope}~{\tt\small  \{oabdelaal, nozbek, jviola, ychen53\}@ucmerced.edu }}%
\thanks{ \small
$^{2}$Adana Alparslan Türkeş Science and Technology University, Faculty of Engineering, Department of Electrical Electronics Engineering, Adana, TÜRKİYE
        \pscircled{\faEnvelope}~{\tt\small nozbek@atu.edu.tr}}%
}


\maketitle

\begin{abstract}
From the perspective of control theory, the gradient descent optimization methods can be regarded as a dynamic system where various control techniques can be designed to enhance the performance of the optimization method. In this paper, we propose a prescribed finite-time convergent gradient flow that uses time-varying gain nonlinear feedback that can drive the states smoothly towards the minimum. This idea is different from the traditional finite-time convergence algorithms that relies on fractional-power or signed gradient as a nonlinear feedback, that is proved to have finite/fixed time convergence satisfying strongly convex or
 the Polyak-Łojasiewicz (PŁ) inequality, where due to its nature, the proposed approach was shown to achieve this property for both strongly convex function, and for those satisfies Polyak-Łojasiewic inequality. Our method is proved to converge in a prescribed finite time via Lyapunov theory. Numerical experiments were presented to illustrate our results. 

\end{abstract}


\section{Introduction}

Optimization plays a crucial role in ensuring efficiency, adaptability, and performance in various applications \cite{Sun2020, GAMBELLA2021_optimization_survey}. In this context, convergence rate and robustness are pivotal considerations for optimization, which are also desirable issues to consider 
for both control theory and machine learning algorithms (e.g.   support vector machines, neural networks, and regression, etc.) where the   convex optimization is widely used \cite{Hu2017_ACC, Lessard2016}.  

The shift in perspective of analyzing the gradient based techniques within a continuous-time framework offers valuable insights in optimization algorithms analysis and synthesis \cite{lin1912control}. A key advantage of continuous-time models is their ability to simulate accelerated approaches, including Nesterov's accelerated gradient descent. The work reported in \cite{Weijie2016} represents a breakthrough in the interpretation of optimization algorithms. Furthermore, continuous analogs of methods such as the Bregman Lagrangian and Polyak's heavy ball provided a model that demonstrate their geometric nature and complex structures \cite{park2024}. From this perspective, the gradient flows emerge from differential equations that capture the dynamics inherent in optimization processes, which constitutes a powerful analytical framework. Furthermore, this perspective allows us to employ Lyapunov designs, which take cues from control techniques such as finite/fixed time control. 
Recently, increasing interest has been seen in the fixed-time convergence design with the capability of ensuring increased speed of convergence, as well as improved transient response.

Gradient-based optimization algorithms rely on asymptotic convergence, which in the majority of applications is insufficient, such as real-time control, online learning, or time-critical decision-making, where solutions need to be accurate and delivered within a timely manner. This encourages researchers to create algorithms with fixed-time and prescribed-time optimization, in which the optimization problem is solved in a predetermined amount of time.

Finite-time and fixed-time optimal optimization methods have emerged to bridge this challenge by ensuring convergence within a finite time period. However, most of the existing finite-time algorithms adopts fractional-power gradient feedback or the signed gradient inspired from the sliding mode control technique, which with aid of Lyapunov constructions are proven to converge in a finite or fixed time for a certain function properties mainly; the function is strongly convex or satisfies the Polyak-Łojasiewicz (PŁ) inequality. Moreover, classical finite-time schemes typically lack explicit control over the exact time of convergence, which is extremely crucial in cases where strict time requirements are necessary, such as robotics trajectory optimization, energy-efficient control, and distributed multi-agent coordination.


Recent studies have highlighted the importance of achieving finite/fixed time convergence \cite{polyakov2012nonlinear, liu2022overview, Song2023_PrTControl_Survey, Polyakov2020}. Some insights into finite-time stability of autonomous systems provided by Bhat and Bernstein \cite{bhat2000finite}. A comprehensive analysis of stability concepts is presented in Polyakov and Fridman \cite{POLYAKOV2014}. Levant and Yu introduced a novel $k$-th order differentiator to address the classical sliding mode control (SMC) problem \cite{Levant2018_SMCkthorder}. With these observations in mind, the finite-time convergence principle is particularly crucial for obtaining correct information within a bounded time frame. By leveraging advancements in the analysis of finite-time control for nonlinear systems, a priori guarantees of finite-time convergence can be provided for optimization algorithms within a deterministic control-theoretic setting. 

Finite-time convergent gradient flows are a class of dynamical systems whose trajectories of optimization variables converge to the optimum in finite time. Unlike asymptotic approaches, finite-time convergence ensures convergence to the solution within a finite time interval, which has some advantages when applications are time-sensitive. The traditional finite-time approaches typically have one limitation of a fundamental flaw: the convergence time is dependent on initial system conditions. This dependence also makes it impossible to predict or establish the exact time of convergence, which, in cases where precision is required about the timing, may not be suitable.

The analysis of finite/fixed-time stability of optimization algorithms has primarily been done for specific problem. In this line, Chen et al. present a finite-time gradient approach by considering the benefits of sliding mode control \cite{ chen2021finite}, wherein the convergence time depends on initial conditions. Then, to make the convergence time robust to initial circumstances, two different designs of fixed-time gradient algorithms are shown. One was created by utilizing the sine function property, which states that the frequency of a sine function determines its convergence time. The other was constructed using the Mittag-Leffler function property, which states that the first positive zero of a Mittag-Leffler function determines its convergence time. In another study, Zhang et al. investigates discrete-time optimization algorithms for finite-time convergent flows \cite{zhang2022firstorderoptimizationinspiredfinitetime}. The technique proposes Euler discretization for rescaled-gradient flow (RGF) and signed-gradient flow (SGF). They analyze convergence guarantees in both deterministic and stochastic settings. From a control application point of view, Su and Chunya proposed a novel finite-time optimization control method for stabilizing switched nonlinear systems with state and control constraints \cite{Su2022_FnT_Optimization}. A switching law is given to ensure finite-time stability across subsystems. The algorithm ensures fast convergence while optimizing system performance and minimizing energy consumption.

Fixed-time stability has been exploited to construct accelerated gradient algorithms for optimization problems in ML of large scale. For instance, Budhraja et al. \cite{Budhraja2021BreakingTC} introduced a gradient-based optimization framework with a fixed-time convergence providing a robustness analysis to noisy gradient, and introducing a discretized model that preserves the fixed time convergence property. On the other hand, Garg et al. \cite{Garg2021_FxTimeFradientFlow} addressed new gradient-flow schemes ensuring fixed-time convergence. The technique, which provides rigorous stability analysis under relaxed convexity conditions, outperformed conventional gradient-based methods in terms of computational efficiency. Similarly, Hu et al. address a modified Projection Neural Network (PNN) to solve nonlinear projection equations with fixed-time convergence \cite{ HU2022_Modified_NN_FxT}. This study verified the existence, uniqueness, and stability of the proposed method using Lyapunov analysis. Toward this goal, strict monotonicity and Lipschitz continuity assumptions are used to ensure stability. Notably, their approach ensures fixed-time convergence. The primary benefits are convergence speed and higher accuracy compared to current PNN techniques. Ozaslan and Jovanovic highlighted a unified framework for modifying exponentially stable optimization algorithms to achieve finite- or fixed-time stability \cite{ozaslan2024exponentialfinitefixedtimestabilityapplications}. The work demonstrated that a simple scaling transformation of the right-hand side of an optimization algorithm can enforce fixed-time stability. The work in \cite{Niroomand2024_FxT_GradientFlow} introduced a new class of gradient flows (GFs) to develop a fixed-time reinforcement learning (RL) algorithm based on a single-network adaptive critic (SNAC). The proposed technique provides a tunable upper bound on convergence time, improving predictability in optimization. Hong et al. developed fixed-time optimization algorithms by utilizing Hessian-inverse-based schemes for solving time-varying convex optimization problems \cite{Hong2023_FxT_TimeVaryingConvex}. They provided a generalizable framework for real-time applications in control and robotics.



This paper introduces a Prescribed Finite-Time Convergent Gradient Flow (PFT-GF) using Time-Varying Feedback, drawing control-inspired feedback concepts into optimization. The key motivation in doing this is the aspiration to leverage dynamical system principles to formally shape the convergence behavior of gradient-based flows. The fundamental idea is to build a gradient flow with a convergence time not only as finite but also as user-prescribed, i.e., the system converges precisely at a pre-specified time point irrespective of the initial conditions. Instead of employing static gain tuning, our method employs a time-varying feedback strategy that modulates the gradient flow trajectory in real time to achieve prescribed finite-time convergence. The approach is inspired from a predefined-time stabilization control, it brings a new appproach of designing optimization algorithms for predefined-guaranteed purposes. The developed approach is particularly valuable in applications requiring fast adaptation and high precision, such as robotics, power systems, and distributed optimization.

The main contributions of this can be summarized as follows: 
\begin{itemize}
\item We introduce a systematic theoretical framework designed by a time-varying feedback mechanism for dynamically adjusting the gradient feedback weight to achieve prescribed finite-time convergent, where due to its nature, it is proven to guarantee this property for both strongly convex function, and for those satisfies Polyak-Łojasiewic inequality.
\item We introduce the sufficient conditions for prescribed finite-time convergence and analyze the stability properties of the resulting system by using the Lyapunov-based convergence proofs.
\item We present a realizable algorithm from the derived gradient flow, which can be applied to a wide range of optimization problems. The algorithm is efficient computationally and easy to implement, and hence it can be used for practical applications.
\item We show various simulation scenarios to demonstrate the practicability and effectiveness of the proposed method. Numerical results verify that the new method achieves exact convergence at the desired time with better performance than the traditional finite-time and asymptotic methods.
\end{itemize}

The remainder of the manuscript is organized as follows: Section II reviews some preliminaries on finite-, fixed-, and prescribed-time stability and convergence methods. Section III provides the mathematical formulation and theoretical analysis of the novel PFT-GF model. Section IV presents simulation results to demonstrate its effectiveness. Finally, Section V gives some concluding remarks and directions for future work.


\section{Preliminaries}

To highlight the characteristics of the proposed gradient flow framework, which ensures convergence within a specified time period, we introduce the necessary definitions and perform a stability analysis in this section.

\begin{definition}
A nonlinear dynamical system is presented as:
\begin{equation}
    \dot{x}(t) = g(x(t), t), 
\end{equation}  
wherein, $g:\mathcal{U} \rightarrow \mathbb{R}^n$ is continuous on an open neighborhood $\mathcal{U} \subseteq \mathbb{R}^n$, and $\quad g(0, t) = 0, \quad x(0) = x_0, \quad t \in \mathbb{R}_{\ge 0}$.

If the system exhibits uniform asymptotic stability, then for every \( x_0 \in \mathcal{U} \) (with \( \mathcal{U} \subseteq \mathbb{R}^n \)) there exists a function \( T_s(x_0) \) such that
\begin{equation}
x(x_0, t) = 0, \quad \forall t \geq T_s(x_0),
\end{equation}
then the system’s equilibrium at the origin is termed \textbf{locally uniformly finite-time stable}. Here, \( T_s(x_0) \) represents the \textbf{settling-time function}.  

\end{definition}

\begin{lemma} \cite{bhat2000finite}  
Consider system (1), where the equilibrium is \textbf{finite-time stable} if there exists a Lyapunov function to be $V(x,t) \in C^1$, $V(x) > 0, \quad \forall x \neq 0, \quad \text{and} \quad V(0) = 0.$ and defined over the domain \( \mathcal{U} \times \mathbb{R}_{\ge0} \), where \( \mathcal{U} \subseteq \mathbb{R}^n \) represents a neighborhood around the origin. Additionally, there exist constants \( \kappa > 0 \) and \( 0 < r < 1 \) such that the function satisfies the following condition:  
\begin{equation}
\dot{V}(x,t) \leq -\kappa V(x,t)^r, \quad \forall x \in \mathcal{U}.
\end{equation}  
Under this condition, for $\forall x_0 \in \mathcal{U}$, the system state reaches the origin in a finite time \( T^* \), which is upper-bounded by:  
\begin{equation}
T^* \leq \frac{V(x_0,t)^{1-r}}{\kappa(1-r)}.
\end{equation}  
If this property extends to the entire state space (i.e., \( \mathcal{U} = \mathbb{R}^n \)), the system exhibits \textbf{global finite-time stability}.  
\end{lemma}

\begin{definition} \cite{POLYAKOV2014, SONG2017_Timevarying_feedback_PrTC}  Consider a system where \( g: \mathcal{U} \times \mathbb{R}_{\geq0} \to \mathbb{R}^n \) represents a nonlinear function that remains locally bounded within an open neighborhood \( \mathcal{U} \subseteq\mathbb{R}^n \) of the origin, and is independent of time in a uniform manner. The state variable is denoted by \( x \in \mathcal{U} \). 
The system is said to exhibit \textbf{fixed-time stability} if there exists a positive constant \( T^* \) such that every trajectory \( x(t, x_0) \) satisfies:  
\begin{itemize}
    \item The system's state converges to equilibrium in a finite duration, implying the existence of a function \( T_s(x_0) \) where
    \begin{equation}
   \hspace*{-8mm}  \lim\limits_{t \to T_s(x_0)} x(t, x_0) = 0, \quad \text{and} \quad x(t, x_0) = 0, \quad \forall t \geq T_s(x_0).
    \end{equation}
    \item The convergence time \( T_s(x_0) \) is upper-bounded by \( T^* \) for all initial conditions \( x_0 \).  
\end{itemize}
\end{definition}


We now provide prescribed-time control, which further extends fixed-time stability by allowing the system designer to set an exact convergence time, making the system reach equilibrium at a user-defined time instant. The use of the monotonically growing function is crucial in the prescribed-time control scheme \cite{SONG2017_Timevarying_feedback_PrTC}. Let us begin by defining the prescribed- time stability.


\begin{definition}\label{DefPT} The origin of the system $\dot{x}(t) = g(x(t), t), $ is known as prescribed/arbitrary time stable if it is fixed time stable, and $\exists \tau_p \in \mathbb{R}_{\geq0}  $, with no dependence on initial conditions. If $T_s(t_0, x_0)=\tau_p$, then the origin is strictly prescribed-time convergent, while $T_s(t_0, x_0) < \tau_p$ indicates weakly prescribed-time convergence.
\end{definition}

Consider a monotonically- incearing time function takes the form:
\begin{equation}\label{MTF}
\mathcal{T}(t-t_0) = \frac{T_p^r}{(T_p + t_0 - t)^r}, \quad t \in [t_0, t_0 + T),
\end{equation}
where $T_p > 0, r \in \mathbb{N}$ and $\mathcal{T}:[t_0,t_0+T_p] \rightarrow [1, \infty]$, note that $\mathcal{T}_1(T_p) = +\infty$, where $T_p$ will be used as the pre-designed time of convergence. Two definitions of prescribed-time stability are reviewed.

\begin{definition} \cite{SONG2017_Timevarying_feedback_PrTC}
A system \( \dot{x} = g(x,t,L) \) where \( x \) and \( L \) are arbitrary-dimensional variables, is said to possess \textbf{Prescribed-Time Input-to-State Stability (PT-ISS)} if there exist class \( \mathcal{KL} \) functions \( \beta \) and \( \gamma \) such that:  
\begin{equation}  
    |x(t)| \leq \beta \left( |x_0|, \mathcal{T}(t - t_0) - 1 \right) + \gamma \left( \sup\limits_{d \in [t_0, t]} |L(d)| \right).
\end{equation}  
Since \( \mathcal{T}(t-t_0) - 1 \) increases monotonically to infinity as \( t \to t_0 + T_p \), a PT-ISS system exhibits global asymptotic stability within time \( T_p \) in the absence of disturbances.
\end{definition}


\begin{definition}\label{def5} \cite{SONG2017_Timevarying_feedback_PrTC}

A system \( \dot{x} = g(x,t,L) \) is said to exhibit \textbf{PT-ISS} if there exist functions \( \beta \) and \( \beta_f \) belonging to class \( \mathcal{KL} \), as well as a function \( \gamma \) from class \( \mathcal{K} \), such that the following conditions hold:

\begin{equation}
|x(t)| \leq \beta_f \left( \beta(|x_0|, t - t_0) + \gamma \left( ||L||_{[t_0,t]} \right), \mathcal{T}(t-t_0) - 1 \right).
\end{equation}
A system that satisfies PT-ISS possesses that its state reaches zero within the prescribed time $T_p$, even in the presence of perturbations.

\end{definition}


\begin{lemma}\label{Lem2} \cite{SONG2017_Timevarying_feedback_PrTC} \newline
Consider the function $\mathcal{T}(t - t_0) = \left( \frac{T}{T + t_0 - t} \right)^{r}$, where $r$ are positive integers. If a continuously differentiable function $V: [t_0, t_0 + T) \to [0,+\infty)$ satisfies
\begin{equation}
    \dot{V}(t) \leq -2\rho \mathcal{T}(t - t_0) V(t) + \frac{\mathcal{T}(t - t_0)}{4\lambda} L(t)^2,
\end{equation}
for positive constants $\rho$ and $\lambda$, then
\begin{equation}
    V(t) \leq e^{-2\rho\int_{t_0}^{t}\mathcal{T}(\tau-t_0)d\tau} V(t_0) + \frac{||L||^2_{[t_0,t]}}{8\rho\lambda}
\end{equation}
with the conditions of lemma 3, if $L(t) \equiv 0$, then  $\lim_{t \to t_0 + T} V(t) = 0$.
\end{lemma}


The existing finite/fixed-time gradient flows use fractional powers of the gradient norms or the sign function inspired from sliding mode control, to guarantee finite/fixed-time convergence. These methods are generally inspired by the q-rescaled continuous gradient flow q-RGF, as well as the q-signed
flow (q-SGF) given as: 

\begin{subequations}
   \begin{equation}
            \dot{x}=F_{q-RGF}(x) = -c \frac{\nabla{f(x)}}{\|\nabla{f(x)}\|^\frac{q-2}{q-1}}
        \end{equation}
             \begin{equation}
            \dot{x}=F_{SGF}(x) = -c \frac{\nabla{f(x)}}{\|\nabla{f(x)}\|^\frac{1}{q-1}}\text{sign}(\nabla{f}(x))
        \end{equation}
	 \label{eq.gf2}
\end{subequations}
\noindent where $c >0, q >1$. Both algorithms are proven to be finite-time convergent, given that the function $f$ is gradient dominated of order $p \in (1,q)$. While for strongly convex functions, both q-RGF and q-SGF is finite-time convergent for any $q \in (2,\infty)$
\cite{romero2020finite}.



However, in this paper we present an alternative to these approaches, a time-varying gradient feedback framework that employs a scaling of the gradient by a function of time that grows monotonically towards the terminal time of convergence in a prescribed finite time.

\section{Design of a PFT-GF with time-varying Gradient feedback}

Consider the following system with a single integrator:

\begin{equation}\label{sysu}
\dot{x} = u 
\end{equation}

\noindent where $x$ and $u$ stand for the control input and system state, respectively. We aim to drive the state to its equilibrium point within a time period $[t_0,T_p]$, where $T_p$ is a predefined time of regulation. The stability analyses for the closed loop are only valid within that time frame. Lets define the following time scaled state feedback control law:
\begin{equation}
    u=-\rho(t) \mathcal{T}(t) x
\end{equation}

Defining a new transformation variable $h=\mathcal{T} x$, the system can be expressed as:
\begin{equation}
\dot{h} = \mathcal{T} u + \dot{\mathcal{T}} x
\end{equation}
where $\mathcal{T}$ is defined as in (\ref{MTF}), and $\dot{\mathcal{T}} = \frac{rT_p^{r}}{(T_p+t_0-t)^{1+r}}$. Then, one can choose a Lyapunov function $V=\frac{1}{2}h^2$.
Upon choosing the time-varying state feedback controller with the gain $\rho = \rho_0 + \frac{r}{T_p} > 0$ \cite{SONG2017_Timevarying_feedback_PrTC}, the derivative of \(V(t,x)\) satisfies the differential inequality:
\begin{equation}
\begin{split}
\dot{V} = h \dot{\mathcal{T}} x + h \mathcal{T} u\\
\le h \mathcal{T}\bigg [ u+\frac{1+r}{T_p}h\bigg]\\
=-\rho_0\mathcal{T}h^2 = -2\rho_0\mathcal{T}V(t).
\end{split}
\end{equation}
Using Lemma \ref{Lem2} , we get:
\begin{equation}
V(t) \leq e^{-2\rho_0 \int_{t_0}^{t}\mathcal{T}(\tau-t_0)d\tau} V(t_0).
\end{equation}

Thus, by Definition \ref{def5}, the trajectory of the system reaches its equilibrium points according to the solution:
\begin{equation}
|x(t)| \leq \mathcal{T}^{-1} e^{- 2\rho_0 \int_{t_0}^{t}\mathcal{T}(\tau-t_0)d\tau} |x(t_0)|. 
\label{PTCsysSol}
\end{equation}

Since $\mathcal{T}^{-1}=1-\frac{t-t_0}{T_p} \rightarrow 0$ as  $t \rightarrow {t_0+T_p}$, then the state is regulated within a prescribed time $T_p$. Note that, the solution (\ref{PTCsysSol}) of the differential equation of the system (\ref{sysu}) converges to its origin as an exponential function of the monotonically increasing time scaling function of the state feedback.

Now, we are tackling the problem of minimizing a function $f: \mathbb{R}^n \rightarrow \mathbb{R} $, in an unconstrained optimization problem setting given as:

\begin{equation}
\begin{aligned}
\min_{x\in  \mathbb{R}^n} \quad & f(x).\\
\end{aligned}
\end{equation}
The cost function given by $f$ is assumed to satisfy:
 The function $f$ has its minimum value
$f^*=f(x^*) =\min f(x) > -\infty$ at $x^* \in \mathbb{R}^{n}$, i.e., $f^* > -\infty$.

Next, we show that the proposed approach for guiding the optimization trajectory to the minimum with a pre-specified time can be achieved for both functions that is strongly convex, and those satisfies the Polyak-Łojasiewicz (PL) inequality. Recall that, strong convexity ensures PL, but not the opposite \cite{karimi2016linear}.
\subsection{The function satisfies the Polyak-Łojasiewicz inequality}

\begin{definition}\label{PL}. Consider a nonempty set $\mathcal{D} \ne \emptyset$, a differentiable function $f : \mathbb{R}^n \rightarrow \mathbb{R} $ is described to be $f \in PŁ_\sigma(\mathcal{D})$ satisfying the Polyak-Łojasiewicz (PŁ) inequality on $\mathcal{D}$, if there exists $r > \min f$ and $\sigma > 0$ such that:
\begin{equation}\label{PLequ}
   [f(x)-\min f]  \le \frac{1}{2\sigma}\|\nabla{f}(x)\|^2  \quad, \forall x \in [\min f < f < r].
\end{equation}
\end{definition}

This inequality indicates that the growth of the gradient is lower bounded by a quadratic function while we are moving far from the minimum function point. The inequality simply implies that each equilibrium point is assumed to be a global minimum, unlike Strong Convexity which implies that there is a unique minimizer, so this inequality can be seen as a weakening of strong convexity.

As a consequence of Definition (\ref{PL}), The function $f \in C_{1,1}^{loc}(\mathbb{R}^n,\mathbb{R})$ has a unique
minimizer $x = x^* $ and $f \in PŁ_\sigma(\mathcal{D})$, as a consequence of gradient dominance, the function $f(x)$ has a a lower bound quadratic growth given as (\cite{karimi2016linear}, Theorem 2):
\begin{equation}
    [f(x)-f^*] \ge \frac{\nu}{2}\|x-x^*\|^2 \quad , \forall x \in \mathbb{R}^n.
\end{equation}

In order to modify the classical gradient flow to guarantee the prescribed finite-time convergence of the optimal point. a simple time-dependent scaling of the original dynamics is introduced, that can be viewed as a time-varying nonlinear feedback of the gradient. A Lyapunov-based scheme is used to prove the prescribed finite-time convergence. 

Consider a GF given by:
    \begin{equation}\label{systv}
        \dot{x}=-k(t) \mathcal{T}(t) \nabla{f}(x(t))
    \end{equation}
where the function  $f$ satisfies assumptions above, $k(t)$ is an adaptive gain, and $\mathcal{T}(t)$ is a monotonically growing function (\ref{MTF}) that plays a crucial role in the prescribed-
time control scheme. 

Choose the positive-definite Lyapunov function $V=[f-f^*]$, the derivative of the function:
\begin{equation}
\begin{split}
\dot{V} = \nabla{f}(x)^T\dot{x}=-k(t)\mathcal{T}(t) \nabla{f}(x)^T \nabla{f}(x)\\=-k(t)\mathcal{T}(t)\|\nabla{f}(x)\|^2
\end{split}
\end{equation}

Since $f \in PŁ_\sigma(\mathcal{D})$, one can apply the inequality condition in (\ref{PL}):
\begin{equation}
\dot{V} \le -2\sigma k(t) \mathcal{T}(t)(f-f^*)=-2\sigma k(t) \mathcal{T}(t)V.
\end{equation}

\noindent where $T_p$ is the prescribed time convergence, $t_0$ is assumed to be the starting time, adopting in the paper to be zero.

With the aid of Lemma (\ref{Lem2})
\begin{equation}
V(t) \leq e^{-2\sigma k(t) \int_{t_0}^{t}\mathcal{T}(\tau-t_0)d\tau} V(t_0)
\end{equation}
and thus,
\begin{equation}
    \lim_{t \rightarrow {t_0+T_p}} V(t)=  0.
\end{equation}
This implies that, regardless of the initial condition, the convergence of the Gradient Flow in (\ref{systv}) can be guaranteed in the prescribed time $T_p$, and thus can be explicitly prespecified.

Now, we show that the approach satisfies this property for a $\mu-$strongly convex function.
\subsection{The function is $\mu-$strongly convex}


\begin{definition}\label{STCX} For a convex function $f(x)$ whose gradient exists, if there exists a scalar $\mu > 0$ such that
\begin{equation}
\langle \nabla f (\chi_1)-\nabla f (\chi_2) , \chi_1 - \chi_2 \rangle \geq \mu\|\chi_1 - \chi_2\|^2,
\end{equation}
for any $\chi_1$ and $\chi_2$ belonging to the definition domain of $f(x)$, then $f(x)$ is said to be $\mu$-strong convex.
\end{definition}

 If the function $f$ is $\mu$-strong convex, then the time scaling approach can drive the states to its minimum point $x^*$ in a prescribed-finite time $T_p$. Consider the Lyapunov function $V =\|x - x^*\|^2$. Taking the time derivative,
\begin{equation}
\begin{split}
\dot{V} = -2k \mathcal{T}(t) {(x - x^*)^T \nabla f(x)}\\
\le -2k \mathcal{T}(t) \|x - x^*\|^2 \\
\leq -2 k \mu \mathcal{T}(t)  V.
\end{split}
\end{equation}
With the aid of Lemma (\ref{Lem2}, $\lim_{t \rightarrow {t_0+T_p}} V(t)=  0$. This implies that, the convergence of the Gradient Flow in (\ref{systv}) can be guaranteed in the prescribed time $T_p$ for  $\mu$-strong convex functions.

\section{Experiments}

\begin{figure} [h]
    \centering
    \includegraphics[width=0.7\linewidth]{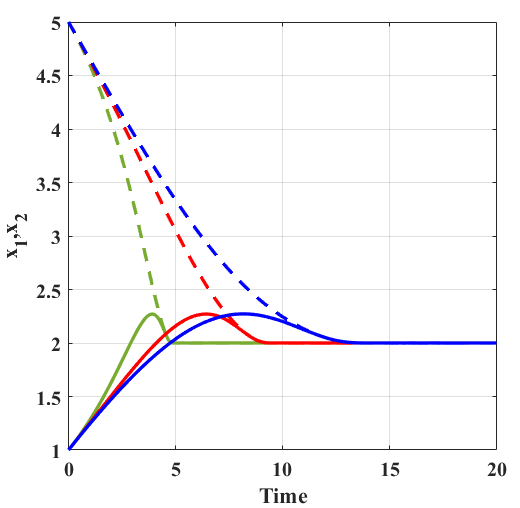}
    \caption{Convergent of the proposed GF for the Trid function ($x_1$ solid, $x_2$ dashed) for different prescribed times $T_p=5,10$ and $15$.}
    \label{fg1l}
\end{figure}

\begin{figure} [h]
    \centering
    \includegraphics[width=0.7\linewidth]{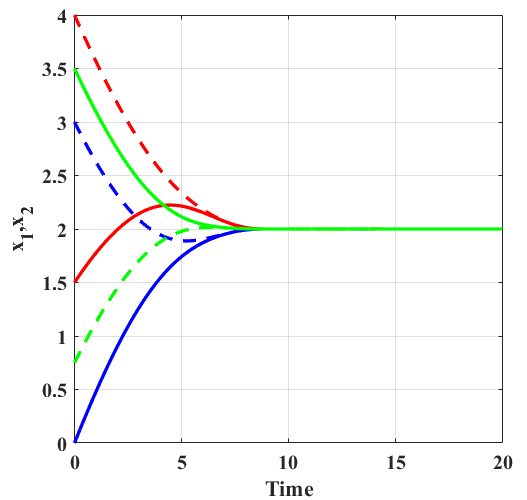}
    \caption{Convergent of the proposed GF for the Trid function from different initial conditions ($x_1$ solid, $x_2$ dashed), $T_p=10,k=0.1$.}
    \label{fg2l}
\end{figure}

\begin{figure} [h]
    \centering
    \includegraphics[width=0.7\linewidth]{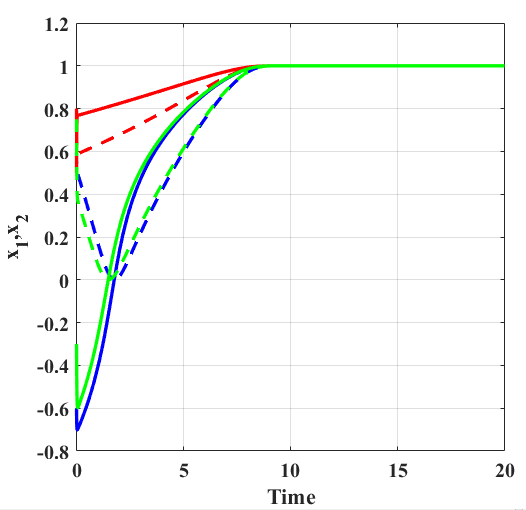}
    \caption{Convergent of the proposed GF for the Rosenbrock function to the global minimum $(1,1)$ from different initial conditions ($x_1$ solid, $x_2$ dashed), $T_p=10,k=0.05$.}
    \label{fg3l}
\end{figure}

To numerically illustrate the prescribed-time convergence of the proposed methods, the convergence behavior of Trid function $f(x) = \sum_{i=1}^{n} (x_i - 1)^2 - \sum_{i=2}^{n} x_i x_{i-1}$ is shown in Figures (\ref{fg1l}, \ref{fg2l}), where Figure (\ref{fg1l}) shows the convergence for three different prescribes time $T_p=5,10$ and $15$. Figure (\ref{fg2l}) illustrate the convergence in a time frame of $T_p=10$ starting from various initial conditions. The Rosenbrock benchmark function $  f(\mathbf{x}) = \sum_{i=1}^{n-1} \left[ 100 (x_{i+1} - x_i^2)^2 + (1 - x_i)^2 \right]$ was evaluated. This non-convex function achieves its global minimum at $(1,1)$. The function is known to be $f\in PL_\sigma([-1, 1] \times [-1, 1])$ with a modulus $\sigma=0.1$. Figure (\ref{fg3l}) shows the convergence of the $x_1,x_2$ towards the global minimum in a prescribed time of $T_p=10$. 

.


\section{Conclusion}

In this paper, we proposed a different Prescribed Finite-Time Convergent Gradient Flow (PFT-GF) with time-varying feedback for enhanced optimization convergence. Unlike traditional finite-time optimization methods with fractional-power gradient feedback, our approach takes a time-dependent scaling function to ensure the system converges to the optimal solution within a given time independent of initial conditions.
The theoretical framework, supported by Lyapunov-based stability analysis, verifies the guaranteed finite-time convergence of the approach. This time-varying control methodology ensures that system behavior can be accurately dictated by the user-defined convergence time without dependence on initial conditions or system parameters. The results suggest potential applications in time-critical areas like real-time control, robotics, and large-scale distributed optimization. Future studies would attempt to extend this methodology to constrained optimization, multi-agent systems and machine learning fields in order to broaden its scope of applications.


\bibliographystyle{ieeetr}
\bibliography{refs}

\end{document}